\documentstyle[12pt]{article}
\catcode`\@=11

\@addtoreset{equation}{section}
\@addtoreset{equation}{subsection}
\def\theequation{\ifnum\value{subsection}>0\relax
\thesubsection.\arabic{equation}\relax
\else\ifnum\value{section}>0\relax
\thesection.\arabic{equation}\relax
\else\arabic{equation}\fi\fi}

\newcommand{\ba}{\begin{eqnarray}}   
\newcommand{\ea}{\end{eqnarray}}     
\newcommand{\be }{\begin{equation}}   
\newcommand{\ee }{\end{equation}}     

\newcommand{\bet}{\beta}
\newcommand{\ga}{\gamma}

\newcommand{\dl}{\delta}

\newcommand{\lam}{\lambda}
\newcommand{\om}{\omega}

\begin{document}
\rightline{  }
\vskip 1cm
\centerline{\Large \bf The Algebraic Structure of the} 
\vskip 0.8cm
\centerline{\Large \bf ${\bf gl(n|m)}$ Color Calogero-Sutherland Models}

\vskip 0.9cm
\centerline{\bf Guo-xing JU$^{1,3}$, Shi-kun WANG$^2$, Ke WU$^1$ }

\vskip0.5cm
\centerline{\sf 
1.Institute of Theoretical Physics}
\centerline{\sf 
Academia Sinica}
\centerline{\sf 
Beijing 100080,China}          
\centerline{\sf 
2.CCAST(World Laboratory),Beijing 100080,China \\
and \\
Institute of Applied Mathematics}
\centerline{\sf  
Academia Sinica}
\centerline{\sf  
Beijing 100080,China}                               
\centerline{\sf   
3.Physics Department}
\centerline{\sf  
Henan Normal University}
\centerline{\sf  
Xinxiang,Henan Province 453002,China}                                                       
\vskip0.5cm
\begin{abstract}
We extend the study on the algebraic structure of the $su(n)$ color Calogero-Sutherland models to the case of $gl(n|m)$ color CS model and show that the generators of the super-Yangian $Y(gl(n|m))$ can  be obtained from two $gl(n|m)$ loop algebras. Also, a super $W_{\infty}$ algebra for the SUSY CS model is constructed.
\end{abstract}
\vskip0.5cm
\newpage

\section{ Introduction}

    The Calogero-Sutherland (CS) models \cite{cal,sat} are one-dimensional particle systems in which all particles interact through inverse square pairwise interactions. The systems are integrable in both classical and quantum cases and have many interesting properties(e.g. the particles have fractional statistics\cite{hal}).\\

    There are various generalizations to the CS models. One of the examples is the $su(n)$ color generalization in which particles have internal degrees of freedom with $n$ possible values \cite{min}. The $su(n)$ extension of the CS models is also integrable and has a new symmetry---Yangian $Y(su(n))$\cite{ber,dri}. For the $su(n)$ CS system, K.Hikami and M. Wadati \cite{hik1,hik2,hik3} have made a distinction between the Calogero-type model and the Sutherland-type model and investigated in detail their integrability and algebraic structure. Hikami and Wadati found that the $W_{\infty}$ algebra, which is the underlying symmetry of the quantum integrable system with long-range interactions, 
unifies the Yangian in the Sutherland-type model and the loop algebras in the Calogero-type model and revealed that Yangian is a subalgebra of the deformed $W_{\infty}$ algebra.\\
    
   Recently, Yangians associated with the simple Lie superalgebras, which are called super-Yangians, have been studied \cite{naz,zha,ahn,ju}. Based on the super-Yangian $Y(gl(n|m))$, a supersymmetric extension of $su(n)$ color CS model is constructed and its integrability has been proved \cite{ahn}. Now the question to be asked is whether the algebraic structure and some properties for $su(n)$ CS model discussed by Hikami, Wadati and others still prevail in this supersymmetric extension. In this paper, we will make a  study to the algebraic structure of the SUSY CS model. \\

The paper is organized as follows. In section 2, the definition of the super-Yangian $Y(gl(n|m))$, the  SUSY CS models, the Lax pairs for these models  and some notations are given. We then study in section 3 the relations between the super-Yangian $Y(gl(n|m))$ for the SUSY Sutherland-type model and two non-commuting loop algebras of which one is for the SUSY Calogero-type model. Moreover, we also derive a super $W_{\infty}$ algebra from the Lax operators for the SUSY Calogero-type model. Finally, some discussions and remarks are made in section 4.

\section{ Supersymmetric Calogero-Sutherland Model and Super-Yangian $Y(gl(n|m))$}

We begin with the definition of the $su(n)$ color Calogero-type and Sutherland-type models. Assume that there are N particles in these systems. Their positions are described by $x_i(i=1,\cdots,N)$ and each particle carries a color degree of freedom labeled by index $a, a=1,\cdots,n$. The Hamiltonians for the $su(n)$ color Calogero-type and Sutherland-type models are respectively defined as \cite{min}
\be
H_c=\frac{1}{2}\sum_{i}(\partial_i)^2+\frac{\lambda}{2} \sum_{i\neq j} (P_{ij}\partial_i\om_{ij}
+\lambda \om_{ij}\om_{ji}),
\ee

\be
H_s=\frac{1}{2}\sum_{i}(x_i\partial_i)^2+\frac{\lambda}{2} \sum_{i\neq j} (P_{ij}+\lambda) \frac{x_i x_j}{(x_i-x_j)(x_j-x_i)},
\ee
where 
\ba
\om_{ij}=\frac{1}{x_i-x_j},&
\theta_{ij}=\frac{x_i}{x_i-x_j},
\ea
and $P_{ij}$ is an permutation operator that interchanges the color degrees of freedom of the $i$-th and $j$-th particles. We denote by $E^{ab}_i, a, b=1,\cdots, n$, the matrices which act as $|a><b|$ on the color degrees of freedom of the $i$th particle and leave the other particles untouched. It is then straightforward to check that
\be
P_{ij}=\sum_{a,b=1}^{n}E_i^{ab}E_{j}^{ba}.
\ee
The matrix unit $E^{ab}_i$ can furnish a vector representation for the Lie algebra $su(n)$.\\

We extend the $su(n)$ color CS models to the $gl(n|m)$ case, that is, the particles in the systems now carry the graded color indices. Let  $V$ be an $n+m$ dimensional $Z_2$ graded vector space and
$\{v^a,a=1,\cdots,n+m\}$ be a homogeneous basis with  grading being defined as follows
\be 
p(a)=\left\{\begin{array}{ll}
           0&\hspace{2cm}a=1,\cdots,n\\
           1&\hspace{2cm}a=n+1,\cdots,n+m.
\end{array}\right.
\ee 
The graded matrix unit $E^{ab}$ is defined as
\begin{equation}
E^{ab}v^c=\delta_{bc}v^a .
\end{equation}
It is clear that we have a vector representation for the Lie superalgebra $gl(n|m)$  given by
\begin{equation}
\left[E^{ab}\right.,\left.E^{cd}\right\}=
\delta_{bc}E^{ad}-(-1)^{\eta}\delta_{ad}E^{cb},
\end{equation}
where $\eta =(p(a)+p(b))(p(c)+p(d))$ and the graded bracket is defined by
\ba
\left[E^{ab},E^{cd}\right\}&
\equiv&E^{ab}E^{cd}-(-1)^\eta
E^{cd}E^{ab}. 
\ea
If we consider N copies of the matrix units $E_i^{ab}(i=1,\cdots,N)$
that act on the i-th space of the tensor product of graded vector spaces $V_1\otimes\cdots\otimes V_N$ with $V_i\cong V$, then we can show that the permutation operator $P_{ij}$ defined as
\be
P_{ij}=\sum_{a,b=1}^{n+m}(-1)^{p(b)}E_i^{ab}E_{j}^{ba}
\ee
exchanges the basis vectors $v_i,v_j$ of $i,j$ spaces and has the following properties
\be
P_{ij}=P_{ji},\hspace{1cm}P_{ij}P_{ij}=1,\hspace{1cm}{\rm and}\hspace{1cm}
P_{ij}P_{jk}=P_{ik}P_{ij}.
\ee
\\

As before, we consider the one-dimensional system with N identical particles. $x_i (i=1, \cdots, N)$ and index $a$ denote respectively the positions and the color degrees of freedom of the particles, but now $a$ is $Z_2$ graded as defined in eq.(2.5) and takes values $a=1, \cdots, n+m$. If we replace $P_{ij}$ in eqs. (2.1) and (2.2) by that defined in eq.(2.9), then we get two Hamiltonians that govern the dynamics of the systems with graded internal degrees of freedom. In this sense, we call the models described by $H_s$ and $H_c$ SUSY Sutherland-type and SUSY Calogero-type models, respectively, similar to what discussed in Ref.\cite{hik1}. \\

The SUSY CS models are integrable. To  prove the integrability of the systems, we may either directly construct the integrals of motion or find a Lax pair. In the former, the system is integrable if the integrals of motion commute among themselves and also commute with the Hamiltonian of the system; while in the latter, the model is integrable if two matrices $L_{ij}$ and $M_{ij}$ with operator entries obey
\be
[H, L_{ij}]=\sum _k (L_{ik}M_{kj}-M_{ik}L_{kj}).
\ee
In Ref.\cite{ahn}, the integrability of the SUSY Sutherland-type model is guaranteed through the construction of the conserved quantities. Here we prove the integrability of SUSY CS models by giving the Lax pairs. For the SUSY Sutherland-type model, the Lax pair is given by
\be
L_{ij}=\delta_{ij}(x_i \frac{\partial}{\partial x_i}+\frac{1}{2})+\lambda(1-\delta_{ij})
\theta_{ij}P_{ij},
\ee
\be
M_{ij}=-\delta_{ij}\lambda\sum_{k\not= i}\theta_{ik}\theta_{ki}P_{ik}+\lambda(1-\delta_{ij})\theta_{ij}\theta_{ji}P_{ij};
\ee
While for the SUSY Calogero-type model, a possible choice is as follows:
\be
L_{ij}=\delta_{ij}\frac{\partial}{\partial x_i}+\lambda(1-\delta_{ij})\om_{ij}P_{ij}.
\ee
\be
M_{ij}=-\delta_{ij}\lambda\sum_{k\not= i}\om_{ik}\om_{ki}P_{ik}+\lambda(1-\delta_{ij})\om_{ij}\om_{ji}P_{ij}.
\ee
$M_{ij}$  in eq.(2.13) and eq.(2.15) satisfies conditions:
\be
\sum_i M_{ij}=\sum_j M_{ij}=0.
\ee
\\

The SUSY Sutherland-type model possesses the super-Yangian $Y(gl(n|m))$\cite{ahn}. Let 
$\{T_p^{ab}\;,p\geq 0\;, a,b=1,\cdots,n+m\}$ be the generators of the super-Yangian $Y(gl(n|m))$, then  $T_p^{ab}$ satisfies the following relations\cite{naz,zha,ahn,ju}
\be
\left[T_s^{ab}\right.,\left.T_{p+1}^{cd}\right\}
-\left[T_{s+1}^{ab}\right.,\left.T_p^{cd}\right\}
=\lambda(-1)^{(p(c)p(a)+p(c)p(b)+p(b)p(a))}\left(T_p^{cb}T_s^{ad}-T_s^{cb}T_p^{ad}\right)
\ee
for $s,p\geq -1$, where $T_{-1}^{ab}\equiv {\lambda}^{-1}{\bf 1}\delta_{ab}$. 
It is easy to check that only $T_0^{ab},T_1^{ab}$ are the basic operators, while  all other $T_p^{ab}(p>1)$ are defined recursively from  $T_0^{ab},T_1^{ab}$ \cite{ju}. Therefore, we need to know the action of operators $T_0^{ab},T_1^{ab}$ on the corresponding configuration space of the system in order to prove that SUSY Sutherland-type model has symmetry $Y(gl(n|m))$. The action of the operators $T_0^{ab},T_1^{ab}$ on the configuration space of  the system is of the following form:
\ba
& &T^{ab}_0=\sum_{i=1}^N E^{ab}_i ,\\
& &T^{ab}_1=\sum_{i,j=1}^N E^{ab}_i L_{ij},
\ea
where $L_{ij}$ is given by eq. (2.12). In addition, we also find that  
\be
T^{ab}_p=\sum_{i,j}E^{ab}_i(L^p)_{ij} \hspace{1.5cm}p\geq 0
\ee
satisfies the relation (2.17). It is straightforward to prove that generators
$T_0^{ab},T_1^{ab}$  are conserved operators for the SUSY Sutherland-type model, that is,
\be
[T_0^{ab}, H_s ]=[T_1^{ab}, H_s ]=0 .
\ee
\\

The SUSY Calogero-type model is a rational limit of the the SUSY Sutherland-type model\cite{ahn}. It only has $gl(n|m)$ loop algebra symmetry which we will discuss in the next section.

\section{ The Algebraic Structure of the SUSY \\Calogero-Sutherland Model}

In this section, we first discuss the defining relations of $Y(gl(n|m))$ in terms of the action of $T_0^{ab}, T_1^{ab}$ on the configuration space of SUSY Sutherland-type model, then we study the loop algebra structure of the SUSY Calogero-type model and its relation with $Y(gl(n|m))$. In the second part of this section, a super $W_{\infty}$ algebra constructed from the Lax operator for the SUSY Calogero-type model is derived.\\

From the section 2, we know that the SUSY Sutherland-type model possesses the super-Yangian $Y(gl(n|m))$ generated by operators $T_0^{ab}, T_1^{ab}$ in eqs. (2.18) and (2.19). Using eq.(2.7), it can be shown that operators $T_0^{ab}, T_1^{ab}$  satisfy the following relations
\ba
  & & [ T_0^{ab} , T_0^{cd} \}  =  \delta_{bc} T_0^{ad} -(-1)^{\eta} \delta_{da}
  T_0^{cb}  \\
  & & [ T_0^{ab} , T_1^{cd} \}  =  \delta_{bc} T_1^{ad} -
  (-1)^{\eta}\delta_{da} T_1^{cb} , \\
  & & [ T_1^{ab} , T_1^{cd} \}  =  \delta_{bc} T_2^{ad} -
  (-1)^{\eta}\delta_{da} T_2^{cb}-\lambda (-1)^{\bet}(T_0^{ad}T_1^{cb}-T_1^{ad}T_0^{cb}),
\ea
where $\bet=p(b)p(c)+p(c)p(d)+p(b)p(d)$ and $T_2^{ab}$ can be written explicitly as
\be
T_2^{ab}=\sum_{i}E_{i}^{ab}(x_i \partial_i +\frac{1}{2})^2 +\lambda \sum_{i\neq j}(E_iE_j)^{ab}(x_i\partial_i\theta_{ij}+\theta_{ij}(x_i\partial_i+x_j\partial_j+1))+\lambda^2\sum_{i\neq j,j\neq k}(E_iE_jE_k)^{ab}\theta_{ij}\theta_{jk}.
\ee
Here we apply the conventional notations, $(E_iE_j)^{ab}=\sum_{c=1}^{n+m}(-1)^{p(c)}E_i^{ac}E_j^{cb}, (E_iE_jE_k)^{ab}=\sum_{c,d=1}^{n+m}(-1)^{p(c)+p(d)}E_i^{ac}E_j^{cd}E_k^{db}$. We can also obtain the following Serre-like relation for operators $T_0^{ab},T_1^{ab}$ 
\ba
& & [ T_0^{ab} , [ T_1^{cd} , T_1^{ef} \} \} -
  [T_1^{ab} , [T_0^{cd} , T_1^{ef} \}\} = \lambda(\delta_{bc}O_{ef}^{ad}-\delta_{de}O_{cf}^{ab} +(-1)^{\delta}\delta_{cf}O_{ed}^{ab}    \nonumber \\
& &+(-1)^{\eta}\delta_{be}O_{af}^{cd} -(-1)^{\eta}\delta_{ad}O_{ef}^{cb} -(-1)^{\ga}\delta_{af}O_{eb}^{cd}),
\ea
with $\delta=(p(c)+p(d))(p(e)+p(f)),\ga=(p(a)+p(b))(p(c)+p(d)+p(e)+p(f))$ and $O_{cd}^{ab}$ being defined as 
\be
O_{cd}^{ab}=-(-1)^{\bet}(T_0^{ad}T_1^{cb}-T_1^{ad}T_0^{cb}).
\ee
The eqs.(3.1)--(3.3) and eq.(3.5) are the defining relations for the super-Yangian $Y(gl(n|m))$.\\

For the SUSY Calogero-type model (2.1), we introduce two operators:
\ba
& &J_0^{ab}=T_0^{ab} ,\\
& &J_1^{ab}=\sum_{i,j}E_i^{ab}I_{ij},
\ea
where
\be
I_{ij}=\delta_{ij}\frac{\partial}{\partial x_i}+\lambda(1-\delta_{ij})\om_{ij}P_{ij}.
\ee
is a Lax operator given in eq.(2.14), here we use another notation $I_{ij}$ for $L_{ij}$. It is easy to check that $J_0^{ab}, J_1^{ab}$ satisfy the following commutation relations
\ba
  & & [ J_0^{ab} , J_1^{cd} \}  =  \delta_{bc} J_1^{ad} -(-1)^{\eta} \delta_{da}
  J_1^{cb} , \\
  & & [ J_1^{ab} , J_1^{cd}  \}  =  \delta_{bc} J_2^{ad} -(-1)^{\eta} 
  \delta_{da} J_2^{cb} , \\
  & & [ J_0^{ab} , [ J_1^{cd} , J_1^{ef} \} \} -
  [ J_1^{ab} ,[ J_0^{cd} , J_1^{ef} \} \}= 0 ,
\ea
where
\be
J_2^{ab}=\sum_{i,j}E_i^{ab}(I^2)_{ij}.
\ee
In general, we have
\ba
& &J_p^{ab}=\sum_{i,j}E_i^{ab}(I^p)_{ij} \hspace{1.5cm} (p\geq 0),\\
& &[ J_s^{ab} , J_p^{cd} \}  =  \delta_{bc} J_{s+p}^{ad} -(-1)^{\eta} \delta_{da}
  J_{s+p}^{cb} .
\ea
The eq.(3.12) is the Serre relation for the $gl(n|m)$ loop algebra. We can show that $J_p^{ab}$ are conserved operators for the SUSY Calogero-type model, that is, the SUSY Calogero-type model has the $gl(n|m)$  loop algebra symmetry. Futhermore, there is another representation of the $gl(n|m)$  loop algebra with the generators $K_p^{ab}$ being defined by
\be
  K_p^{ab} = \sum_{i=1}^N E_i^{ab} x_i^p  \hspace{1.5cm} (p\geq 0),
\ee
and having the relations
\ba
&  & [ K_s^{ab} , K_p^{cd} \}  =  \delta_{bc} K_{s+p}^{ad} -(-1)^{\eta}
    \delta_{da}  K_{s+p}^{cb} ,\\
&  &[ K_0^{ab} , [ K_1^{cd} , K_1^{ef} \} \} -
  [ K_1^{ab} ,[ K_0^{cd} , K_1^{ef} \} \}= 0 .
\ea
The $K_0^{ab},K_1^{ab}$ are the two basic operators with the following forms, respectively
\ba
&  &K_0^{ab}=J_0^{ab} ,\\
&  &K_1^{ab}=\sum_iE_i^{ab}x_i.
\ea
\\

We consider now the algebras constructed from the elements $\{J_0^{ab},J_1^{ab},K_1^{ab}\}$. Just as in the non-graded case\cite{hik1,ber1}, the super-Yangian $Y(gl(n|m))$ appears by unifying two $gl(n|m)$ loop algebras $\{J_0^{ab},J_1^{ab}\}$ and $\{K_0^{ab},K_1^{ab}\}$, that is to say, $T_1^{ab}$ can be obtained from the commutators between two sets of generators $\{J_0^{ab},J_1^{ab}\}$ and $\{K_0^{ab},K_1^{ab}\}$:
\ba
&  &[ J_1^{ab} , K_1^{cd} \}+ [ K_1^{ab} , J_1^{cd} \} -\lambda[ J_0^{ab} , J_0^{cd} \}+   \nonumber \\
&  &+\lambda [ J_0^{ab} , (J_0J_0)^{cd} \}= 2
  (\delta_{bc} T_1^{ad} - (-1)^{\eta}\delta_{da} T_1^{cb}) .
\ea
In this sense, the algebra generated by $\{J_0^{ab},J_1^{ab},K_1^{ab}\}$ is a larger algebra including both super-Yangian $Y(gl(n|m))$ and two loop algebras.
In addition, we have two extra Serre-like relations:
\ba
 & & [ J_0^{ab} , [ J_1^{cd} + T_1^{cd} , J_1^{ef} + T_1^{ef} \}\}
   - [ J_1^{ab} + T_1^{ab} , [ J_0^{cd} , J_1^{ef} + T_1^{ef} \}\}
   = \lambda(\delta_{bc}M_{ef}^{ad}-\delta_{de}M_{cf}^{ab} \nonumber \\
& &+(-1)^{\delta}\delta_{cf}M_{ed}^{ab}    
   +(-1)^{\eta}\delta_{be}M_{af}^{cd} 
   -(-1)^{\eta}\delta_{ad}M_{ef}^{cb} 
   -(-1)^{\ga}\delta_{af}M_{eb}^{cd}),\\
& & [ J_0^{ab} , [ K_1^{cd} + T_1^{cd} , K_1^{ef} + T_1^{ef} \}\}
   - [ K_1^{ab} + T_1^{ab} , [ J_0^{cd} , K_1^{ef} + T_1^{ef} \}\}
   = \lambda(\delta_{bc}N_{ef}^{ad}-\delta_{de}N_{cf}^{ab}  \nonumber \\
& & +(-1)^{\delta}\delta_{cf}N_{ed}^{ab}    
    +(-1)^{\eta}\delta_{be}N_{af}^{cd} 
    -(-1)^{\eta}\delta_{ad}N_{ef}^{cb} 
    -(-1)^{\ga}\delta_{af}N_{eb}^{cd}), 
\ea
where
\ba
& &M_{cd}^{ab}=(-1)^{\bet}[\sum_{i\neq j}E_i^{ad}E_j^{cb}((x_i+1)\frac{\partial}{\partial x_i}- (x_j+1)\frac{\partial}{\partial x_j})
+\lambda \sum_{ijk} \prime\frac{x_i+1}{x_i-x_j}(E_iE_j)^{ad}E_k^{cb} \nonumber \\
& &-\lambda \sum_{ijk}\prime\frac{x_j+1}{x_j-x_k}E_i^{ad}(E_jE_k)^{cb}]
+\lambda \sum_{i\neq j}\frac{x_i+x_j+2}{x_i-x_j}E_i^{ab}E_j^{cd} \\
& &N_{cd}^{ab}=(-1)^{\bet}[\sum_{i\neq j}E_i^{ad}E_j^{cb}(\frac{\partial}{\partial x_i}-\frac{\partial}{\partial x_j}+x_i-x_j)
+\lambda \sum_{ijk}\prime\frac{1}{x_i-x_j}(E_iE_j)^{ad}E_k^{cb} \nonumber \\
& &-\lambda \sum_{ijk}\prime\frac{1}{x_j-x_k}E_i^{ad}(E_jE_k)^{cb}]
+\lambda \sum_{i\neq j}\frac{2}{x_i-x_j}E_i^{ab}E_j^{cd} 
\ea
and $\sum \prime$ means that any two summation indices do not coincide. The Serre-like relations (3.5),(3.22) and (3.23) may be written in a more compact form. To do this, we define the operator $Q_1^{ab}(x,y)$ to be of the following form:
\be
  Q_1^{ab} (x,y) \equiv T_1^{ab} + x \, J_1^{ab} + y \, K_1^{ab},
\ee
where $x$ and $y$ are two  complex number. Using the operator $Q_1^{ab}(x,y)$, we get the following relation
\ba
& &[ J_0^{ab} , [ Q_1^{cd}(x,y) , Q_1^{ef}(x,y) \} \} -
  [Q_1^{ab}(x,y) , [J_0^{cd} , Q_1^{ef}(x,y) \}\} =  \nonumber \\
& &=\lambda(\delta_{bc}P_{ef}^{ad}(x,y) 
-\delta_{de}P_{cf}^{ab}(x,y) +(-1)^{\delta}\delta_{cf}P_{ed}^{ab}(x,y)       
+(-1)^{\eta}\delta_{be}P_{af}^{cd}(x,y)  \nonumber \\
& &-(-1)^{\eta}\delta_{ad}P_{ef}^{cb}(x,y) -(-1)^{\ga}\delta_{af}P_{eb}^{cd}(x,y)),
\ea
where
\ba
& &P_{cd}^{ab}(x,y)=O_{cd}^{ab}+x\{(-1)^{\bet}[\sum_{i\neq j}E_i^{ad}E_j^{cb}(\frac{\partial}{\partial x_i}-\frac{\partial}{\partial x_j})
\nonumber \\
& &+\lambda \sum_{ijk}\prime\frac{1}{x_i-x_j}(E_iE_j)^{ad}E_k^{cb} 
-\lambda \sum_{ijk}\prime\frac{1}{x_j-x_k}E_i^{ad}(E_jE_k)^{cb}]
+ \sum_{i\neq j}\frac{2}{x_i-x_j}E_i^{ab}E_j^{cd}\}   \nonumber \\
& &+y(-1)^{\bet}\sum_{i\neq j}E_i^{ad}E_j^{cb}(x_i-x_j).
\ea
From eq.(3.27), we see that generators $Q_1^{ab}(x,y)$  also form a representation of the super-Yangian $Y(gl(n|m))$ for any $x$ and $y$. Therefore, we can get a family of super-Yangian subalgebras from generators $Q_1^{ab}(x,y)$, all these subalgebras have a common deformation parameter $\lambda$.\\

  Using the generators $J_p^{ab}$  and $K_p^{ab}$, we can also construct another algebra. We introduce two set of operators 
as 
\be 
{\cal W}_p^{(s)}=\frac{1}{2(p+s)}[\sum_i x_i^2,{\cal W}_{p+2}^{(s-1)}],
\ee

\be 
{\cal Q}_p^{(s)ab}=\frac{1}{2(p+s)}[\sum_i x_i^2,{\cal Q}_{p+2}^{(s-1)ab}],s=2,3,\cdots
\ee
where ${\cal W}_p^{(1)}=J_p=\sum_{ij}(I^p)_{ij}, {\cal Q}_p^{(1)ab}=J_p^{ab}$. ${\cal W}_p^{(s)}$ and ${\cal Q}_p^{(s)ab}$ can be considered to have conformal spin $s$. We can write ${\cal W}_p^{(s)}$ and ${\cal Q}_p^{(s)ab}$ explicitly as
\be 
{\cal W}_p^{(s)}=\frac{1}{2^{s-1}(p+s)_{s-1}}[\underbrace{\sum_i x_i^2,[\cdots,[\sum_i x_i^2}_{s-1},J_{p+2s-2}]\cdots ]]=\sum_j(-x_j)^{s-1}(\partial _{x_j})^{p+s-1}+\cdots ,
\ee
\be 
{\cal Q}_p^{(s)ab}=\frac{1}{2^{s-1}(p+s)_{s-1}}[\underbrace{\sum_i x_i^2,[\cdots,[\sum_i x_i^2}_{s-1},J_{p+2s-2}^{ab}]\cdots ]]=\sum_j E_j^{ab}(-x_j)^{s-1}(\partial _{x_j})^{p+s-1}+\cdots ,
\ee
where $(p)_s$ denotes
\be
(p)_s=p(p+1)\cdots (p+s-1),
\ee
and $(\cdots)$ represents the lower order terms of $\partial_{x_j}$ and is not necessarily zero in the limit $\lam\rightarrow 0$. Similar to the proofs by Hikami and Wadati \cite{hik1}, we can show that ${\cal W}_p^{(s)}$ and ${\cal Q}_p^{(s)ab}$ satisfy the following commutation relations:
\be
[{\cal W}_p^{(s)},{\cal W}_q^{(s^\prime)}]=((s-1)q-(s^\prime -1)p){\cal W}_{p+q}^{(s+s\prime -2)}+\cdots,
\ee

\be
[{\cal W}_p^{(s)},{\cal Q}_q^{(s^\prime)ab}]=((s-1)q-(s^\prime -1)p){\cal Q}_{p+q}^{(s+s\prime -2)}+\cdots,
\ee

\be
[{\cal Q}_p^{(s)ab},{\cal Q}_q^{(s^\prime)cd}\}=(\dl _{bc}{\cal Q}_{p+q}^{(s+s\prime -1)ad}-(-1)^\eta\dl _{ad}{\cal Q}_{p+q}^{(s+s\prime -1)cb}) +\cdots,
\ee
where $(\cdots)$ includes lower spin operators. Eq.(3.34) shows that the operators ${\cal W}_p^{(s)}$ constitute the quantum $W_\infty$ algebra as in the CS models\cite{hik1}. However, eqs.(3.35) and (3.36) indicate that algebra generated by operators ${\cal W}_p^{(s)}$ and ${\cal Q}_p^{(s)ab}$ is a kind of super  $W_{\infty}$ algebra with color\cite{berg}.\\

In  the case of $\lambda\longrightarrow 0$, if we define a operator
which is of the following form
\ba
  {\cal Q }^{(s)ab}_{p}=\sum_{i=1}^N E^{ab}_i x_i^{s-1}
  ( \partial_{ x_i} )^{p+s-1} ,
\ea
that is we only consider the first term of ${\cal Q}_p^{(s)ab}$ in eq.(3.32) up to a constant factor, then it is easy to show that ${\cal Q} ^{(s)ab}_{p}$ satisfy the commutation relations:
\ba
  [ {\cal  Q }^{(s)ab}_{p} , {\cal Q} ^{(s')cd}_{q} \} & = &
  \delta_{bc} \cdot \sum_{k=0}^{p+s-1}
\frac{ (p+s-1)!(s'-1)!}{k!(p+s-k-1)!(s'-k-1)!}
  \, {\cal Q} _{p+q}^{(s+s'-1-k)ad} \nonumber \\
  & & \ -(-1)^{\eta} \delta_{da}  \cdot \sum_{k=0}^{q+s'-1}
\frac{ (q+s'-1)!(s-1)!}{k!(q+s'-k-1)!(s-k-1)!}
  \,{\cal Q} _{p+q}^{(s+s'-1-k)cb}  .  \nonumber \\
\ea
It can be seen that the algebra generated by ${\cal Q} ^{(s)ab}_{p}$ in eq.(3.37) is a kind of super  $W_{\infty}$ algebra without central terms\cite{berg}.\\

\section{ Remarks and Discussions}

In this paper, we have discussed the algebraic structure of the SUSY CS models. We know that SUSY CS models have super-Yangian $Y(gl(n|m))$ symmetry, and we give the  realization of $Y(gl(n|m))$ and its relations with two $gl(n|m)$ loop algebras . Also, a super  $W_{\infty}$ algebra realized in terms of generators ${\cal W}_p^{(s)}$ and ${\cal Q}_p^{(s)ab}$ is shown.
Compared with the CS models, there also exists a lager symmetry generated by two loop algebras for the SUSY CS models. But, there are three points to be remarked. Firstly, the operator $L_{ij}$ (eq.(2.14)) is different from the one for the Sutherland-type model in Ref.\cite{hik1,ber1}, so that the right hand sides of the Serre-like relations (3.22) and (3.23) are nonzero. Secondly, in the right hand side of eq.(3.27), there is a $x,y$ dependence which is absent in  the non-graded CS models\cite{hik1,ber1}. Thirdly, the $W_\infty$ algebra is graded with respect to the color indices $a,b$.\\

Because of the above similarity between SUSY CS models and the CS models, we can, in the same way as done in Refs.\cite{hik1,hik2}, discuss the problem of energy spectrum and the symmetry of the system confined in a harmonic potential. But, we don't here make further study on these problems.

\section*{Acknowledgment}
The first of the authors  thanks  Wei-Zhong Zhao for useful discussions,and also acknowledges  K.Hikami for providing  his publications.

\end{document}